\documentclass[article]{asl}
\title[A Model Theory for the Potential Infinite]{A Model Theory for the Potential Infinite}
\author{Eberl, Matthias}
\revauthor{Eberl, Matthias}
\email{matthias.eberl@mail.de}
\address{Mathematisches Institut, LMU,\\
Theresienstr. 39, D-80333 M\"unchen, Germany}

\newtheorem{Theorem}{Theorem}[section]
\newtheorem{theorem}[Theorem]{Theorem}
\newtheorem{proposition}[Theorem]{Proposition}
\newtheorem{lemma}[Theorem]{Lemma}

\newtheorem{corollary}[Theorem]{Corollary}
\theoremstyle{definition}
\newtheorem{definition}[Theorem]{Definition}

\newtheorem{example}[Theorem]{Example}
\newtheorem{remark}[Theorem]{Remark}

\usepackage{graphicx}
\usepackage{proof}

\usepackage{graphicx}
\def\HHH{\mathcal{H}}
\def\III{\mathcal{I}}
\def\JJJ{\mathcal{J}}
\def\KKK{\mathcal{K}}
\def\LLL{\mathcal{L}}
\def\MMM{\mathcal{M}}

\def\PPP{\mathcal{P}}

\def\SSS{\mathcal{S}}
\def\TTT{\mathcal{T}}

\def\VVV{\mathcal{V}}

\def\phi{\varphi}
\def\rho{\varrho}
\def\theta{\vartheta}
\def\vecb#1{\boldsymbol{#1}}
\def\eps{\,\epsilon\,}
\def\epsi{\,\epsilon_{(i,i)}\,}
\def\imp{\Rightarrow}

\def\iffdef{\, {:\!\iff}\,}

\def\embb{\hookrightarrow}
\def\emb#1#2{emb_{#1}^{#2}}

\def\up{\uparrow\!}
\def\val#1{\lbrack\!\lbrack#1\rbrack\!\rbrack}
\def\valp#1#2{\lbrack\!\lbrack{#1}\rbrack\!\rbrack^{#2}}
\def\pot#1{\PPP(#1)}
\def\potfin#1{\PPP_{fin}(#1)}

\def\nat{\mathbb{N}}
\def\real{\mathbb{R}}
\def\nil{()}
\def\iset#1{\III_{#1}(\vecb{a}:C)}
\def\isetp#1#2{\III_{#1}(#2)}

\def\fil{\mathfrak{D}}

\def\valid{\models}
\def\validu#1#2{\models_m#1[#2]}
\def\validun{\models_m}
\def\validb#1#2{\models_{\ll}#1[#2]}

\def\validbn{\models_{\ll}}
\def\validv#1#2{\models#1[#2]}

\def\stdmod{\mathfrak{M}^{Tar}}
\def\allapp{\mathfrak{M}^{ind}}

\def\tdr#1#2{#1 : #2}
\def\tdm#1#2{#1 \mid #2}
\def\tdf#1#2#3{#1 : #2 \to #3}
\def\tdt#1#2#3{#2 \mid #1 : #3}

\def\y#1{\mathsf{#1}}
\def\yx#1{\mathsf{x_{#1}}}
\def\yt#1{\mathsf{t_{#1}}}
\def\yequi{\longleftrightarrow}
\def\yA{\mathsf{\Phi}}
\def\yB{\mathsf{\Psi}}
\def\yBpow{\mathsf{\Psi_{pow}}}
\def\yBpair{\mathsf{\Psi_{pair}}}

\begin{document}
\begin{abstract}
We present the model theoretic concepts that allow mathematics to be developed with the notion of the potential infinite instead of the actual infinite. The potential infinite is understood as a dynamic notion, being an indefinitely extensible finite. The main adoption is the interpretation of the universal quantifier, which has an implicit reflection principle. Each universal quantification refers to an indefinitely large, but finite set. The quantified sets may increase, so after a reference by quantification, a further reference typically uses a larger, still finite set. We present the concepts for classical first-order logic and show that these dynamic models are sound and complete with respect to the usual inference rules. Moreover, a finite set of formulas requires a finite part of the increasing model for a correct interpretation.
\end{abstract}

\maketitle

\section{Introduction}
\label{intro}

The aim of this paper is to present a model theory based on a potentialist's viewpoint, i.e., infinity is understood as a potential infinite. An introduction of this approach has been presented in \cite{eberl2021Logica}. In short, the usual formal language, in this paper first-order predicate logic, is directly interpreted in this model, whereby the domain of a universal quantifier is an indefinitely large, but finite set. A philosophical discussion of this concept has already been expounded by Shaughan Lavine in the Section ``The Finite Mathematics of Indefinitely Large Size'' in his book \cite{lavine2009understanding}, based on Jan Mycielski's work about locally finite theories \cite{Mycielski1986}.

\subsection{Motivation and Related Work.}

In \cite{shapiro2006all} Shapiro and Wright considered the potential infinite as an indefinite extensible concept. Therein they state that ``If a `collection' is not a set, then it is nothing, has no size at all, and so can't be {`too big'}'' and ask: ``The question, simply, is whether it is ever appropriate or intelligible to speak of all of the items that fall under a given indefinitely extensible concept''. The authors come to the conclusion that there is no satisfying solution how to read such a quantification and refer to \emph{reflection principles} as a possible answer. The interpretation that we present has this implicit reflection principle, with the characteristic that the use of the phrase ``for all elements $\dots$'' may immediately change the state of the indefinitely extensible collection.

Let us shortly recap Dummett's understanding\footnote{Whereas Dummett concludes that statements quantifying over an indefinitely extensible concept do not follow the laws of classical logic, our model theoretic approach does not require any modifications of axioms or inference rules.} of this notion in \cite{dummett1994mathematics}: ``An indefinitely extensible concept is one such that, if we can form a definite conception of a totality all of whose members fall under that concept, we can, by reference to that totality, characterize a larger totality of all whose members fall under it.'' The ordinal numbers and sets are a typical example, but also the natural numbers form such an indefinitely extensible concept. If we refer to ``all sets'', this creates or reveals\footnote{We want to emphasize that our terminology does not assume any philosophical view such as Platonism, formalism, intuitionism or intentionalism (the latter was Jan Mycielski's view which he explained in \cite{mycielski1989meaning}).} a new set and thus the totality of all sets has changed. But already the concept ``number of numbers'' is indefinitely extensible in this sense (first there is no number, creating $0$, thus there is one number, creating $1$ and so on). The forthcoming formalization is based on the following line of thoughts:
\begin{enumerate}
\item The infinite is understood as a potential infinite, which is a form of \emph{ontological finitism} without any actual infinite sets.
\item The potential infinite is regarded as being \emph{indefinitely extensible}, which is a \emph{dynamic} concept. Infinite sets are exhausted by procedures.
\item It nevertheless allows \emph{indefinitely large finite} stages as a substitute for completed infinities.
\item These indefinitely large states are the basis for an interpretation with an \emph{implicit reflection principle}.
\end{enumerate}

This reflection principle as part of the interpretation avoids situations as in set theory with unrestricted quantification in the interpretation, which is justified afterward, e.g.~by the L\'{e}vy-Montague's reflection principle or by a hierarchy of Grothendieck universes.

Marcin Mostowski uses potential infinite sets $(\nat_i)_{i \in \nat}$ with $\nat_i := \{0, \dots, i-1\}$ as a basis of his model theory, but his approach is less dynamic and the axioms of Peano arithmetic must be adopted in a way that there exists a greatest number (see~\cite{mostowski2016truth}, \cite{mostowski2003representing}). Carlson \cite{carlson2003ranked} uses a ``ranked'' model with reference to Mycielski's work. Recently Linnebo and Shapiro \cite{linnebo2019actual} suggested to formalize the potential infinite using a modal reading of this notion. They rely on an analysis of the potential infinite worked out by Niebergall \cite{niebergall2014assumptions}. 

All of the mentioned approaches use some translation, restriction or adoptions of the original axioms. Our approach is new, to our knowledge, insofar as it formalizes the notion of an indefinitely large finite as a relation, being part of the interpretation. As a consequence, no manipulation of the axioms is necessary.

\subsection{Outline of the Idea.}
\label{outlsec}

For the purpose of explanation, start from a Tarskian model $(\MMM,\valid)$ with an actual infinite carrier set. In a first step exchange the underlying universe by a potential infinite set. This is formally an increasing family (a direct system) $\MMM_\III := (\MMM_i)_{i \in \III}$ with finite sets $\MMM_i$ and directed index set $\III$. The index set represents the stages, necessary to see a set as being dynamic, which may be read figuratively as time. The preorder $\leq$ on $\III$ expresses that stage $i$ is preceding stage $i'$ whenever $i \leq i'$. The relations on $\MMM$, and similarly functions, are replaced by families of relations with states $R_C \subseteq \MMM_C := \MMM_{i_0} \times \dots \times \MMM_{i_{n-1}}$ for $C = (i_0, \dots, i_{n-1}) \in \III^n$. Whenever $C \leq C'$, then $R_C$ is the restriction of $R_{C'}$ to $\MMM_C$. A single state $\MMM_i$ is in general not a (Tarskian) model. For instance, the standard model for Peano arithmetic has an underlying carrier $(\nat_i)_{i \in \nat}$ with $\nat_i := \{0, \dots, i-1\}$, but none of the sets $\nat_i$ is closed under the successor operation.

The key part of the interpretation is a notion of being \emph{indefinitely large}, given as a relation $C \ll i$ between a list $C = (i_0, \dots, i_{n-1})$ of indices, called \emph{context}, and the indefinitely large index $i$. The state at index $i$ is a substitute for an absolute infinitely large state. The set $\{i \in \III \mid C \ll i\}$ defines an \emph{indefinitely large region} \emph{inside} the system, replacing the single actual infinite state \emph{outside} of $\MMM_\III$. It is enough to test quantified propositions on elements in $\MMM_i$ for such an indefinitely large index $i$ whether these propositions are true or not. New elements in $\MMM_{i'}$ with $i' \geq i$  do not change the truth values of the propositions. 

The interpretation $\valid$ of a formula $\yA$ is then replaced by $\validbn$, with relation $\ll$ as an additional parameter. The variable assignment $\vecb{a}$ is taken from $\MMM_C$, leading to $\validb{\yA}{\vecb{a}:C}$. The basic difference to usual interpretations is the reading of the universal quantifier: $\validb{\forall \y{x} \yA}{\vecb{a}:C}$ holds by definition iff (i.e.,~if and only if)
\[
\validb{\yA}{\vecb{a}b:Ci} \text{ holds for all elements } b \in \MMM_i \text{ for some } i \gg C,
\]
whereby $\vecb{a}b$ denotes the list $\vecb{a}$, extended by $b$ (and $b$ is assigned to the variable $\y{x}$). The reason for this definition is that the use of a potential infinite set does not permit us to speak about all elements in \emph{all} stages $\MMM_i$, because it constitutes an increasing or ``open'' collection and the locution ``all'' has no fixed meaning. 

As a consequence, quantifiers in a formula refer to different stages of an increasing domain. In a formula, say $\forall \yx{0} \, \exists \yx{1} \yA$, the stage to which $\yx{1}$ refers is typically larger than the one to which $\yx{0}$ refers\footnote{Hilary Putnam already used a similar idea in \cite{putnam1967mathematics} with concrete graph models of Zermelo set theory. His aim was to show that an understanding of ``mathematics as a modal logic'' is equivalent to ``mathematics as set theory''. Putnam  proposes an interpretation of a single formula in increasing graph models, without reference to a maximal model, in order to translate a set theoretic reading of a statement into a modal logic one.}. We show in Section \ref{soundcomp2} that the usual classical deduction rules are nevertheless sound and complete w.r.t.~the interpretation in these models.

It is important to note that this form of potentialism is an \emph{ontological} one, not an epistemological one. In particular, Hilbert's Program is concerned with the justification of classical mathematics by finitary reasoning. Constructive mathematics uses finite algorithmic procedures. We do not claim that there is an effective procedure to determine an indefinitely large finite state. We only assume that if an existential assertion is true, then there must be an instance for which the statement is true\footnote{Basically, we refer to the meta-level notion and use the equivalence of the following statements (for a property $P$ on natural numbers): $\exists n \in \nat \ P(n)$, and $\exists n \in \nat_i \ P(n)$ for some $i \in \nat$. This equivalence is used for the interpretation of Section \ref{intinreflp}. Moreover, we assume that existential and universal quantification are reducible to each other. On meta-level we use classical reasoning, on object-level we investigate classical logic. This corresponds to a ``liberal potentialist'' in \cite{linnebo2019actual}, since we do not have a notion of \emph{knowledge at stage $i \in \III$}, e.g.~knowledge of whether $\forall n \in \nat \ P(n)$ holds or whether there is some counterexample. Nevertheless, our approach can be applied to Kripke models and intuitionistic logic as well, see Section \ref{simpadopt}.}. Even if the existential assertion is unbounded, nonetheless, this does not require completed infinite sets: If an existential assertion is true, it makes sense to go step by step until one finds the witness in finitely many steps and then to stop there. Of course, the truth of the statement is not decidable by this procedure, but this is epistemological issue.

Compared to Kripke models, the set of nodes $\KKK$ in a Kripke model is usually not directed. When these models are used in order to show intuitionistic invalid propositions, then they often use different nodes without upper bound. Even more important is the fact that properties and relations in Kripke models might grow although the set of objects is the same. In a model presented here, if a property $P$ does not hold at some stage $i$ for $a$, it will not hold at a later state $i' \geq i$ on $a$ either. Finally, a Kripke model uses a quantification over all (possibly infinitely many) later states $k' \geq k$ to interpret $\forall \y{x} \yA$ at node $k$. By that, the quantification may use infinite sets of objects.

\subsection{A Consequent Finitistic Reading.}
\label{reflsec}

Sets on meta-level, such as the index set $\III$ or relation $\ll$, are in most situations infinite, so the question arises, whether actual infinite sets are still used. Probably the presented concepts cannot be developed without them.

Let us first assume that the infinite sets on meta-level are actual infinite sets. Then the investigation shows how to develop potential infinite models and how to relate them to actual infinite Tarskian models. Such a comparison obviously requires the notion of an actual infinite set, at least the carrier set $\MMM$ of a Tarskian model is actual infinite (if it is infinite). It is natural then to assume that meta-level sets, such as the index set $\III$, are actual infinite as well. This leads to a comparison as formulated in Lemma \ref{multlem}.

Nevertheless, Theorem \ref{finiterestthm} states that one never needs the whole, infinite model $\MMM_\III$ to interpret expressions. So if we do not compare potential infinite models to actual infinite models, we can avoid actual infinity on meta-level as well. There is however an unavoidable circularity, because a consistent presentation requires the understanding of the interpretation from the very beginning. Compared to other approaches, e.g.~\cite{linnebo2019actual}, there is no other language or concept as a ``rock bottom'', instead, the new interpretation affects its own presentation. A way to deal with this situation is that in a first reading, one may regard infinite on meta-level as actual infinite. If the idea, that an infinite set is always an indefinitely extensible system with indefinitely large intermediate states, is clear, a second reading is possible. 

In this second reading one should from the very beginning view infinite sets as dynamic concepts with a context dependent extension. The context is given by the currently used expressions, objects, relations and states of these. The context on meta-level is left implicit and the paper shows that it is always possible to make the context explicit, justifying the assumption that such a context always exists. Applied to the meta-level, one would have to explicitly formulate the background-model and to formalize the language in which the concepts have been developed. As a consequence, the locution $i \in \III$, for instance, is a syntactical expression which has an interpretation as an element $i$ and an increasing predicate $(\III_j)_{j \in \JJJ}$, such that $i \in \III$ means $i \in \III_j$ for a sufficiently large index $j$. Theorem \ref{finiterestthm} then states that the background-model is itself finite\footnote{But we cannot assume that all our sets are fixed finite sets, in particular, we cannot use the property that an index set has a largest element. The index set $\nat_\nat$ has a finite size at any moment of a specific reference, but the set could increase and a further reference uses a larger (still finite) set.} (at any time), but adequate for all the expressions used to develop the model theory and to apply it to specific models of interest. This consequent finitistic approach does not compare the potential infinite models with actual infinite Tarskian models anymore, but dynamic models with intermediate and indefinitely large states of this dynamic model.

Although it is possible to apply the results to first-order ZFC set theory, this is not sufficient for a consequent treatment of mathematics with a potential infinite. A set in ZFC set theory is a single object $a$ that represents a ``real'' set through the interpretation of the membership relation, which is in this model theory a potential infinite set. However, a potential infinite set is primarily not a single entity, but a family of states of this object. Higher-order logic introduces indefinitely extensible objects as objects of a higher type. These are different to single objects. For instance, a Cauchy sequence is an indefinitely extensible sequence of higher type whereas a real number is a single base type object.

\subsection{Structure of the Paper}
\label{structpapsec}

In Section \ref{mult} we develop the replacement of a static carrier set by a system $\MMM_\III$ and relations by families of relations. A $\Sigma$-structure ($\Sigma$ a signature) becomes a $\Sigma^\III$-structure, whereby $\Sigma^\III$ is a refinement of $\Sigma$. The system $\MMM_\III$ reflects the potentialist view on a set (and likewise on a relation) as indefinitely increasing. We also introduce a notion to express that an index set has ``almost all'' indices and extend this to contexts $C = (i_0, \dots, i_{n-1})$. We call these sets \emph{indefinitely large} and they reflect the cardinal aspect of the potential infinite.

Section \ref{indefappsec} introduces the core concepts. First we introduce relation $\ll$, mentioned in Section \ref{outlsec}, which corresponds to the ordinal aspect of infinity. For an infinite structure, $\ll$ is itself indefinitely large. Before we define the interpretation, we introduce an intermediate notion between syntax and semantics, that of an \emph{state declaration} $\tdt{\y{t}}{C}{i}$ (with $i \in \III$) and $\tdm{C}{\yA}$ for terms $\y{t}$ and formulas $\yA$. It uses relation $\ll$ for the quantifier and guarantees that the necessary approximation instances of the functions and relations are available for the interpretation. We define the interpretation only for expressions (i.e., terms and formulas) having such a state declaration, but show that for an indefinitely large model all expressions have such a declaration.

Section \ref{intrefl} presents the main results. To ensure that the interpretation is sound, additional requirements on relation $\ll$ are necessary. Roughly, if $i \gg C$ holds, then the set $\MMM_i$ must contain all witnesses of existential quantified formulas, given that the variable assignment is taken from $\MMM_C$. The construction is similar as in the proof of the L\"owenheim-Skolem theorem. This is handled in Section \ref{defmeans}. In Section \ref{adeqint} we show amongst other topics that the interpretation is independent of the chosen state declaration.

The main result is Corollary \ref{soundcompl2} and Theorem \ref{finiterestthm}: The interpretation with reflection principle is sound and complete and it requires for a finite set of formulas $\TTT$ only a finite substructure $\MMM_\JJJ$ (with $\JJJ$ a finite index set) to interpret all formulas in $\TTT$ correctly --- Mycielski proved this in \cite{Mycielski1986} by translating formulas and using a common Tarskian semantics. In Section \ref{extexsec} we carry out the developed concepts by use of set theory as an example. Section \ref{further} concludes the paper with further remarks, mainly possible modifications.

\section{Indefinitely Extensible Structures}
\label{mult}

This section describes the first step towards a model theory with potential infinite sets, the use of a system instead of the static carrier set of a Tarskian model.

\subsection{Systems.}
\label{dynset}

The mathematical formalization of an \emph{indefinitely extensible structure} is based on a system\footnote{The family $\MMM_\III$ is a special case of a direct system. When extending the approach to higher-order logic, also inverse systems and more general notions of a system are required.}, having several stages. This is a family of sets 
\[
\MMM_\III := (\MMM_i)_{i \in \III},
\]
such that each $\MMM_i$ is finite. The set $\III$ of indices or stages is a non-empty directed set with a preorder $\leq$ such that $\MMM_i \subseteq \MMM_{i'}$ holds for $i \leq i'$. We use the abbreviation $\up i := \{i' \in \III \mid i' \geq i\}$. In order to be in line with the common definition of a model, at least one set $\MMM_i$ must be non-empty.

The main example is $\nat_\nat = (\nat_i)_{i \in \nat}$, whereby the index set $\nat$ is equipped with the usual order $\leq$ and $\nat_i$ denotes $\{0, \dots, i-1\}$. A further example is $\PPP_\nat = (\PPP_i)_{i \in \nat}$ with $\PPP_i = \pot{\nat_i}$. Here $\pot{\nat_i}$ refers to the powerset of $\nat_i$. Moreover, each finite set $\MMM$ is such a family in a trivial way: The index set is any singleton set, say $\III = \{\ast\}$, and $\MMM_\ast := \MMM$. The interesting situation is when $\III$ is unbounded, but we do not exclude the case that $\III$ is a fixed finite set. In that case a greatest element $j \in \III$ and a comprehensive set $\MMM_{j}$ exists.

A list of indices $C = (i_0, \dots, i_{n-1}) \in \III^n$ is called \emph{context}, or more specifically \emph{state context}. The empty context is denoted as $\nil$ and $Ci$ is the result of adding the index $i$ to an existing context $C$. $\MMM_C$ stands for $\MMM_{i_0} \times \dots \times \MMM_{i_{n-1}}$ and $\up C$ for $\up i_0 \times \dots \times \up i_{n-1}$. Moreover, $C' \leq C$ is defined pointwise for two contexts of the same length.

\subsection{Indefinitely Many Elements.}
\label{indefmanysec}

Many of the subsequent concepts need a notion of ``infinitely many'' indices. Cofinality of a set $\JJJ \subseteq \III$ is a weak notion of infinity, e.g., the intersection of two cofinal sets can be empty. What is required is a stronger notion, one which roughly states that finally all elements are in that set. Moreover, this notion must apply to contexts.

Consider a set $\HHH \subseteq \III^{n} \times \III$ of contexts as a $n$-ary multivalued function $C \mapsto \{i \in \III \mid Ci \in \HHH\}$ that provides possible extensions of a context $C \in \III^n$. We are interested in those $\HHH$ that have sufficiently many values $i$ for sufficiently many arguments $C \in \III^n$. This idea is formalized in the next definition by sets $\fil_{n+1}$. Therein \emph{up-set} means a non-empty, upward closed set. 

\begin{definition}
\label{indefmanydef}
The sets $\fil_{n} \subseteq \pot{\III^{n}}$ are defined recursively as follows:
\begin{align*}
\HHH \in \fil_0 \iffdef& \HHH = \{\nil\}. \\
\HHH \in \fil_1 \iffdef& \HHH \text{ contains an up-set, i.e.,} \up i \subseteq \HHH \text{ for some } i \in \III. \\
\HHH \in \fil_{n+1} \iffdef& d_n(\HHH) \in \fil_n \text{ with } d_n(\HHH) := \{C \in \III^n \mid C^\HHH \in \fil_1\} \\
&\text{and } C^\HHH := \{i \in \III \mid Ci \in \HHH\}.
\end{align*} 

A family $\HHH_\nat := (\HHH_n)_{n \in \nat}$ with $\HHH_n \subseteq \III^n$ is \emph{indefinitely large} or has \emph{indefinitely many} contexts iff $\HHH_n \in \fil_n$ for all $n \in \nat$.
\end{definition}

For a family $\HHH_\nat = (\HHH_n)_{n \in \nat}$ we shortly write $C^\HHH$ for $C^{\HHH_{n+1}}$ if $C \in \III^n$. The next proposition is essential, but its proof is straightforward.

\begin{proposition}
\label{intersectlem}
For all $n \in \nat$, the set $\fil_n$ is a (proper) filter. Each set $\HHH \in \fil_n$ is cofinal and $\up C \in \fil_n$ holds for all $C \in \III^n$.
\end{proposition}

Each set $\HHH_n \subseteq \III^n$ can be extended to a family $\HHH_\nat$ in a natural way:

\begin{definition}
\label{gendef}
Let $\HHH_n \subseteq \III^n$, then there is a family $\HHH_\nat$ \emph{generated by} $\HHH_n$: For $m < n$ let $(i_0, \dots i_{m-1}) \in \HHH_m$ iff there are indices $i_{m}, \dots, i_{n-1}$ such that $(i_0, \dots, i_{n-1}) \in \HHH_n$. For $m > n$ let $\HHH_{m}$ be the set $\HHH_{n} \times \III^{m-n}$.
\end{definition}

\begin{lemma}
\label{extlem}
If $\HHH_n \in \fil_n$, then the family $\HHH_\nat$ generated by $\HHH_n$ is indefinitely large. Moreover, 
\begin{equation}
\label{shorteneq}
Ci \in \HHH_{m+1} \text{ implies } C \in \HHH_m \text{ for all } m \in \nat.
\end{equation}
\end{lemma}

\begin{proof}
We show $\HHH_m \in \fil_m$ for $m \leq n$ inductively on $n-m$ and for $m \geq n$ inductively on $m$: The case $m = n$ holds by assumption. For $m < n$ we have $d_m(\HHH_{m+1}) \subseteq \HHH_m$ since $C^\HHH \in \fil_1$ implies $Ci \in \HHH_{m+1}$ for some $i \in \III$. Moreover, $d_m(\HHH_{m+1}) \in \fil_m$ follows from $\HHH_{m+1} \in \fil_{m+1}$ (induction hypothesis), hence $\HHH_{m} \in \fil_{m}$ by the filter property.

Assume $m \geq n$ and we have to show $d_m(\HHH_{m+1}) \in \fil_{m}$. We have $C^\HHH = \III \in \fil_1$ for all $C \in \HHH_{m}$, consequently $\HHH_{m} \subseteq d_m(\HHH_{m+1})$. The claim then follows from the induction hypothesis $\HHH_{m} \in \fil_{m}$. Property (\ref{shorteneq}) is obvious.
\end{proof}

Operations and relations on families of sets are defined pointwise, for instance, the intersection $\HHH_\nat \cap \HHH_\nat'$ is $(\HHH_n \cap \HHH'_n)_{n \in \nat}$.

\begin{lemma}
\label{indefcor}
If $\HHH_\nat$ and $\HHH_\nat'$ are indefinitely large, then $\HHH_\nat \cap \HHH_\nat'$ is indefinitely large, too. If $\HHH_\nat$ and $\HHH_\nat'$ satisfy Property (\ref{shorteneq}), so does  $\HHH_\nat \cap \HHH_\nat'$.
\end{lemma}

\subsection{The Structure on Systems.}
\label{unimultsec}

In order to make the system $\MMM_\III$ a first-order model, we need to define functions and relations on it. To simplify the presentation we use relations only. 

\begin{definition}
\label{dircomplemma}
A family $R_\HHH := (R_C)_{C \in \HHH}$ of $n$-ary relations $R_C \subseteq \MMM_C$ with $\HHH \subseteq \III^n$ is called \emph{compatible} iff $R_C(\vecb{a}) \iff R_{C'}(\vecb{a})$ holds for all $\vecb{a} \in \MMM_C \cap \MMM_{C'}$. A relation on a system $\MMM_\III$ is a compatible family $R_\HHH$ with a non-empty index set $\HHH$. $R_\HHH$ is \emph{indefinitely large} iff $\HHH \in \fil_n$.
\end{definition}

Compatibility implies that $R_C$ is simply the restriction of $R_{C'}$ to $\MMM_C$ if $C \leq C'$. The requirement $\HHH \not= \emptyset$ in this definition is quite weak. However, we do not require the stronger condition of being indefinitely large in the general definition, since finite structures do not necessarily satisfy them (cf.~Section \ref{smallsub}).

The notion of a signature $\Sigma$ in the common Tarskian semantics is refined to a signature $\Sigma^\III$. It consists of assignments $\tdr{\y{R}}{C}$ for relation symbols $\y{R}$ in $\Sigma$, if $C$ is of length $arity(\y{R})$, such that each relation symbol $\y{R}$ has at least one assignment $\tdr{\y{R}}{C}$. If not mentioned otherwise, the signature $\Sigma^\III$ simply contains $\y{R} : C$ for \emph{all} state contexts $C$ of length $arity(\y{R})$.

\begin{definition}
\label{signdef}
Given a signature $\Sigma^\III$, based on $\Sigma$. A \emph{$\Sigma^\III$-structure} is a family $\MMM_\III$ (as introduced in Section \ref{dynset}) and a map assigning to each assignment $\tdr{\y{R}}{C}$ an instance $R_{C} \subseteq \MMM_C$ of a relation $R_\HHH$ (with $C \in \HHH$ and $\HHH \subseteq \III^n$). We call $\Sigma^\III$ indefinitely large\footnote{It may be the case that $\III$ is infinite but $\bigcup_{i \in \III} \MMM_i$ is nevertheless finite. Then the signature $\Sigma^\III$ of a structure $\MMM_\III$ can be indefinitely large by this definition, although the whole structure is finite. This does not lead to problems, however, we can also ignore these kind of structures because, if $\MMM := \bigcup_{i \in \III} \MMM_i$ is finite, they can be replaced by a trivial structure with $\III = \{\ast\}$ and $\MMM_* = \MMM$.} iff $\{C \in \III^n \mid \tdr{\y{R}}{C}\} \in \fil_n$ for each n-ary relation symbol $\y{R}$ in $\Sigma^\III$.
\end{definition}

To each $\Sigma$-structure $\MMM$ there are associated $\Sigma^\III$-structures $\MMM_\III$ and vice versa:

\begin{lemma}
\label{multlem}
Given a $\Sigma$-structure $\MMM$ and a cover $\MMM = \bigcup_{i \in \III} \MMM_i$. If the sets $\MMM_i$ are finite and $\III$ is directed, then $\MMM_\III$ with $i \leq i' \iffdef \MMM_i \subseteq \MMM_{i'}$ defines a $\Sigma^\III$-structure, whereby $\tdr{\y{R}}{C}$ is part of $\Sigma^\III$ for all contexts $C$. The relations $R_C$ are the restrictions of $R$ to $\MMM_C$, hence $\HHH = \III^n$ holds for all $n$-ary relations $R_\HHH$.

Conversely, given a $\Sigma^\III$-structure $\MMM_\III$ with relations $R_\HHH$, then the union of all elements $\MMM := \bigcup_{i \in \III} \MMM_i$ together with relations $R := \bigcup_{C \in \HHH} R_C$ defines a $\Sigma$-structure. 
\end{lemma}

The simplest way to construct a $\Sigma^\III$-structure from a $\Sigma$-structure is by choosing an index set $\III$ which is isomorphic to $\potfin{\MMM}$, the set of all finite subsets of $\MMM$, together with relation $\subseteq$. The example $\nat_\nat$ shows that the index set need not be the whole collection $\potfin{\MMM}$, but it can be a well-ordered part of it. There are however situations in which this is not possible (or requires at least the well-ordering theorem)\footnote{For instance, the set $\potfin{\real}$ is an uncountable collection without a well-ordered subset having $\real$ as its union.}.

\section{State Declarations}
\label{indefappsec}

In this section we give the prerequisites that are necessary to define for a formula $\yA$ the interpretation $\validb{\yA}{\vecb{a}:C}$. A peculiarity of this interpretation is that it can be applied only to formulas $\yA$ having an \emph{state declaration} $\tdm{C}{\yA}$. Intuitively, $\tdm{C}{\yA}$ says that the formula $\yA$ has a meaning relative to context $C$. 

The language $\LLL$ that we consider is that of a first-order predicate logic of signature $\Sigma$. We assume that there is a fixed sequence $\yx{0},\yx{1},\dots$ of variables, the constant $\bot$, the primitive connectives $\to$, $\land$, $\lor$ and quantifiers $\forall$, $\exists$. We apply the usual abbreviations, for instance, $\neg \yA$ for $\yA \to \bot$ and $\top$ for $\neg \bot$. We use relation symbols only (function symbols and terms are treated shortly in Section \ref{simpadopt}). 

To each expression we explicitly add the list of free variables, which is always a list of the form $(\yx{0},\dots,\yx{n-1})$, whereby quantifiers bind the last\footnote{This kind of convention is sometimes called ``inverse (or dual) de Bruijn notation'' or ``de Bruijn levels''. The order $\yx{0}, \yx{1}, \dots$ corresponds to the order of bound variables in a formula which allows a straightforward mapping of variables to positions in a context and in a variable assignment $\vecb{a} = (a_0, \dots, a_{n-1})$. That is, $a_k$ is assigned to $\yx{k}$ in a formula $\yA(\yx{0},\dots,\yx{n-1})$. We especially use this advantage when we present an extended example in Section \ref{extexsec}.} variable of the list. We will write as usual $\yA(\yx{0},\dots,\yx{n-1})$ to indicate this and sometimes we call these expression ``$n$-ary''. We often omit the subscript of a variable, simply writing $\y{x}$.

\subsection{The Indefinitely Large Relation.}
\label{relinf}

The core concept in order to interpret a formula in a potential infinite structure $\MMM_\III$ is the notion of an ``indefinitely large finite'' or a ``relative infinite'', given by the notion of a $\ll$-relation.

\begin{definition}
\label{lldef}
A family $\ll \, := (\ll_n)_{n \in \nat}$ with $\ll_n \subseteq \III^n$ is named \emph{$\ll$-relation} if it satisfies Condition (\ref{shorteneq}), i.e., $Ci \in \, \ll_{n+1}$ implies $C \in \, \ll_n$ for all $n \in \nat$. An element $C \in \, \ll_n$ in a $\ll$-relation is a \emph{$\ll$-context} (of length $n$).
\end{definition}

We use infix notation, i.e., $Ci \in \ll_{n+1}$ is written as $C \ll i$, including $\nil \ll i$ for $i \in \, \ll_1$. Then $C^\ll$ is the set $\{i \in \III \mid C \ll i\}$ and we sometimes write $i \gg C$ instead of $C \ll i$. Note that a $\ll$-context $C = (i_0, \dots, i_{n-1})$ satisfies by definition $(i_0, \dots, i_{k-1}) \ll i_k$ for all $0 \leq k < n$. 

The set $C^\ll$ corresponds to an \emph{indefinitely large region} and if $C^\ll = \ \up h$, then we call $h$ the \emph{horizon}, see the figure below. The intuition is that the number of objects and means at the current stage, given by the context $C$, is finite, and an indefinitely large (or sufficiently large) index $i \gg C$ is beyond the scope that we can overlook from there.

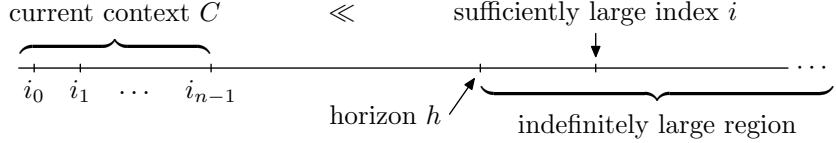
\begin{figure}[ht]
\centering
\includegraphics[trim = -9mm 0mm -3mm -4mm, width=1.0\textwidth]{I.1}
\label{fig1}
\end{figure}

So far we have not defined what exactly it means that an index $i$ is indefinitely large relative to a context $C$. This will be done in Section \ref{sufflarge}. For infinite structures $\MMM_\III$ we will later assume that $\ll$ is indefinitely large in the sense of definition \ref{indefmanydef}. A trivial example of an indefinitely large $\ll$-relation is $\III^* = (\III^n)_{n \in \nat}$. We already have, and we later will call mathematical concepts ``indefinitely large'' if the set of involved states is in the filter $\fil_n$. These are:
\begin{enumerate}
\item Families $\HHH_\nat$ iff all $\HHH_n \in \fil_n$, Def.~\ref{indefmanydef}.
\item $n$-ary relations $R_\HHH$ iff $\HHH \in \fil_n$, Def.~\ref{dircomplemma}.
\item Signatures $\Sigma^\III$ iff $\{C \in \III^n \mid \tdr{\y{R}}{C}\} \in \fil_n$ for all $R:C$ in $\Sigma^\III$, Def.~\ref{signdef}.
\item State declarations $\tdm{C}{\yA}$ iff $\{C \in \III^n \mid \, \tdm{C}{\yA}\} \in \fil_n$, Def.~\ref{typedecldef}.
\item $\Sigma^\III$-models $(\MMM_\III, \validbn)$ iff $\Sigma^\III$ and $\ll$ are indefinitely large, Def.~\ref{termdef}.
\end{enumerate}

One of the task is to show that the notion of being indefinitely large transfers from one concept to the next.

\subsection{Defining the State Declarations.}
\label{typedecltermform}

Assume the model is using the extensible structure $\nat_\nat$ of natural numbers and let $\yA(\yx{0}, \yx{1})$ be the formula $\yx{0} + 1 \leq \yx{1}$. Here we use a function symbol $+$ that we have not introduced yet, but we look at it shortly in Section \ref{simpadopt}. In order to interpret the formula we have to select states $i_0 \in \III$ and $i_1 \in \III$ for the variables $\yx{0}$ and $\yx{1}$. The variable assignment $(a_0,a_1)$ will then be taken from $\nat_{i_0} \times \nat_{i_1}$. To interpret the constant $1$ requires at least $\nat_2$, and to interpret $\yx{0} + 1$ requires at least $\nat_{i_0 + 1}$. The interpretation of formula $\yA$ then uses an instance $\leq_{(i_0 + 1,i_1)} \subseteq \nat_{i_0 + 1} \times \nat_{i_1}$ for the relation symbol $\leq$.

In order to make sure that we have the instances of (functions and) relations that are large enough for an interpretation, we introduce a binary relation between a context $C \in \III^n$ and a formula $\yA(\yx{0}, \dots, \yx{n-1})$, written as $\tdm{C}{\yA}$. For terms $\y{t}$, in particular variables, we introduce a ternary relation $\tdt{\y{t}}{C}{i}$ between a context $C \in \III^n$, a term $\y{t}$ and an index $i \in \III$.

\begin{definition}
\label{typedecldef}
Given a signature $\Sigma^\III$ and a $\ll$-relation $\ll$. The \emph{state declarations} $\tdt{\yx{k}}{C}{j}$ and $\tdm{C}{\yA}$ are relations between a variable $\yx{k}$, a $\ll$-context $C$ of length $n$ and an index $j \in \III$, resp.~between a formula $\yA(\yx{0}, \dots, \yx{n-1})$ and $\ll$-context $C$ of length $n$. State declarations are defined recursively on $\yA$. 
\begin{align*}
\tdt{\yx{k}}{C}{j} \ \iffdef& j \geq i_k \ \text{ for } C = (i_0, \dots,i_{n-1}),\\
\tdm{C}{\y{R}\yt{0} \dots \yt{m-1}} \ \iffdef& \tdt{\yt{0}}{C}{j_0}, \dots, \tdt{\yt{m-1}}{C}{j_{m-1}} \text{ and } \\
& \tdr{\y{R}}{(j_0,\dots,j_{m-1})} \ \text{ for some } j_0, \dots, j_{m-1},\\
\tdm{C}{\yA \to \yB} \ \iffdef& \tdm{C}{\yA} \text{ and } \tdm{C}{\yB} \text{ (and similar for $\land$, $\lor$)}, \\
\tdm{C}{\forall \y{x} \yB} \ \iffdef \tdm{C}{\exists \y{x} \yB} \ \iffdef& \tdm{Ci}{\yB}\text{ holds for some } i \in \III \text{ with } i \gg C.
\end{align*}

The relation $\tdm{C}{\bot}$ holds for all $\ll$-contexts $C$. An expression is called \emph{approximable} (in $\Sigma^\III$ with $\ll$) if there is some state declaration for it. Likewise, a set of expressions is called approximable if all expressions therein are approximable. We call the set of declaration for a formula $\yA$ \emph{indefinitely large} iff $\{C \in \III^n \mid \, \tdm{C}{\yA}\} \in \fil_n$. 
\end{definition}

Remind that by Definition \ref{lldef}, if $Ci$ is a $\ll$-context, then $C$ is a $\ll$-context, too. If the language has no function symbols, the terms $\yt{k}$ in Definition \ref{typedecldef} are all variables. Expressions usually have several state declarations, and for the same declaration $\tdm{C}{\yA}$ there may also be several ways how it has been derived, i.e., there could be different declarations for subformulas.

One might compare a state declaration with a common type declaration $\tdt{\y{t}}{\Gamma}{\sigma}$ of a term $\y{t}$ (having type $\sigma$ in the type context $\Gamma$). A type declaration adds further information to terms in order to rule out meaningless terms and to interpret them properly. Similarly, a state declaration $\tdt{\y{t}}{C}{i}$ adds the necessary state information that is necessary to interpret the terms dynamically (and similarly the formulas). In both cases the declarations are constraints to yield meaningful expressions only, as well as to define their meaning. The next lemma roughly states that in an infinite structure all formulas are approximable with infinitely many contexts.

\begin{lemma}
\label{apxtypdecllem2}
Assume the signature $\Sigma^\III$ and the $\ll$-relation are indefinitely large. Then the set of declarations $\tdm{C}{\yA}$ is indefinitely large for all formulas $\yA$. In particular, each formula is approximable.
\end{lemma}

\begin{proof}
By induction on $\yA$. For atomic formulas we show that $\tdm{C}{\y{R}\yt{0} \dots \yt{m-1}}$ holds for each $\ll$-context $C$ of length $n$. Firstly, the set $\{C' \in \III^m \mid \y{R} : C'\}$ is in $\fil_m$ by assumption. Secondly, the set of all $(j_0, \dots, j_{m-1})$ with $\tdt{\yt{0}}{C}{j_0}$, $\dots$, $\tdt{\yt{m-1}}{C}{j_{m-1}}$ is an up-set and hence in $\fil_m$. Since their intersection is in $\fil_m$, and thus non-empty, we have a context $C' = (j_0, \dots, j_{m-1})$ with $\y{R} : C'$ and $\tdt{\yt{k}}{C}{j_k}$. This shows $\tdm{C}{\y{R}\yt{0} \dots \yt{m-1}}$.

The claim for the connectives follows directly from the filter property (Proposition \ref{intersectlem}). For the quantifier we have by definition 
\begin{align*}
\{C \in \III^n \mid \, \tdm{C}{\forall \y{x} \yB}\} &= \{C \in \III^n \mid \exists i \in \III \text{ with } \tdm{Ci}{\yB} \text{ and } C \ll i\} \\
&\supseteq \{C \in \III^n \mid \{i \in C^\ll \mid \, \tdm{Ci}{\yB}\} \in \fil_1\} \\
& = \{C \in \III^n \mid \{i \in \III \mid \, \tdm{Ci}{\yB}\} \in \fil_1\} \cap \{C \in \III^n \mid C^\ll \in \fil_1\}.
\end{align*}

In order to show $\{C \in \III^n \mid \, \tdm{C}{\forall \y{x} \yB}\} \in \fil_n$ it suffices to show that both sets on the r.h.s.~are in $\fil_n$. By induction hypothesis $\{Ci \in \III^{n+1} \mid \, \tdm{Ci}{\yB}\} \in \fil_{n+1}$ and hence $\{C \in \III^n \mid \{i \in \III \mid \, \tdm{Ci}{\yB}\} \in \fil_1\} \in \fil_n$. The set $\{C \in \III^n \mid C^\ll \in \fil_1\}$ is in $\fil_n$ since $\ll_{n+1} \in \fil_{n+1}$ holds by assumption.
\end{proof}

The next definition is used in Section \ref{smallsub}. Its main application is for an infinite set $\III$ with a finite subset $\JJJ \subseteq \III$.

\begin{definition}
\label{sufflargedef}
Let $\TTT$ be a set of formulas, all approximable in $\Sigma^\III$ with $\ll$. Then a directed subset of indices $\JJJ \subseteq \III$ is a \emph{possible restriction for $\TTT$} iff all formulas in $\TTT$ are approximable in $\Sigma^\JJJ$ with the restriction of $\ll$ to $\JJJ$.
\end{definition}

From this definition it follows immediately that any set $\JJJ' \subseteq \III$ with $\JJJ \subseteq \JJJ'$ is again a possible restriction for $\TTT$.

\begin{lemma}
\label{finiterestlem}
Given a finite set of formulas $\TTT$, approximable in $\Sigma^\III$ with $\ll$, as well as a finite set of indices $\III_0 \subseteq \III$. Then there is a finite set $\JJJ$ with $\III_0 \subseteq \JJJ \subseteq \III$ such that $\JJJ$ is a possible restriction for $\TTT$.
\end{lemma}

\begin{proof}
All formulas in $\TTT$ are approximable in $\Sigma^\III$, so let $\TTT^+$ contain for each formula $\yA \in \TTT$ an (arbitrarily chosen) state declaration $\tdm{C}{\yA}$. For a state declaration $\tdm{C}{\yA}$ there is a finite set of involved indices, defined recursively on Definition \ref{typedecldef}:
\begin{align*}
\JJJ_{\tdt{\yx{k}}{C}{j}} &:= \{i_0, \dots, i_{n-1},j\},  && C = (i_0, \dots, i_{n-1}), \text{ and } j \geq i_k,\\
\JJJ_{\tdm{C}{\bot}} &:= \{i_0, \dots, i_{n-1}\},  && C = (i_0, \dots, i_{n-1}),\\
\JJJ_{\tdm{C}{\y{R} \yt{0} \dots \yt{m-1}}} &:= \bigcup_{0 \leq k < m} \JJJ_{\tdt{\yt{k}}{C}{j_k}},  && \text{with }\tdr{\y{R}}{(j_0, \dots, j_{m-1})},\\
\JJJ_{\tdm{C}{\yA \to \yB}} &:= \JJJ_{\tdm{C}{\yA}} \cup \JJJ_{\tdm{C}{\yB}}, &&\text{similar for } \yA \land \yB \text{ and } \yA \lor \yB,  \\
\JJJ_{\tdm{C}{\forall \y{x} \yA}} := \JJJ_{\tdm{C}{\exists \y{x} \yA}} &:= \JJJ_{\tdm{Ci}{\yA}}, && \text{with } i \gg C.
\end{align*}

The required index set $\JJJ$ is $\bigcup_{\tdm{C}{\yA} \in \TTT^+} \JJJ_{\tdm{C}{\yA}} \cup \III_0 \cup \{j\}$ with $j \in \III$ an upper bound of this set. The sole role of index $j$ is to make $\JJJ$ a directed set. The way the indices have been selected guarantees that the state declaration of $\tdm{C}{\yA} \in \TTT^+$ is also a state declaration in $\Sigma^\JJJ$.
\end{proof}

In some situations it is possible to select a single $\ll$-context of length $n$ such that each $n$-ary formula is approximable with this single context:

\begin{example}
\label{simpex}
For each relation symbol $\y{R} \in \Sigma$ let the signature $\Sigma^\III$ contain $\tdr{\y{R}}{C}$ for all contexts $C$ of the length $arity(\y{R})$. Moreover, $\Sigma$ shall have no function symbols. Assume that $\ll$ is such that we can define functions $\iota_n$, selecting a specific index from $C^\ll$ ($n$ the length of $C$) as follows:
\[
i_0 := \iota_0(\nil^\ll), \ i_1 := \iota_1((i_0)^\ll), \ i_2 := \iota_2((i_0,i_1)^\ll), \dots
\]

This is for instance possible if $\ll$ is indefinitely large and $\ll_n \, \subseteq d_n(\ll_{n+1})$ holds, i.e., $C^\ll \in \fil_1$ for all $\ll$-contexts $C$. The latter may not be the case in general. Then each context $(i_0,\dots, i_{n-1})$ is a $\ll$-context and
\[
\tdm{(i_0,\dots, i_{n-1})}{\yA(\yx{0}, \dots, \yx{n-1})}
\]
holds for all formulas $\yA$. For any set of formulas $\TTT$, the set $\{i_0, i_1, i_2 \dots\}$ is a possible restriction for $\TTT$. If $\TTT$ is a finite set such that all variables (free and bound) are within $\yx{0}, \dots, \yx{n-1}$, then the set $\{i_0, \dots, i_{n-1}\}$ is already a possible restriction for $\TTT$.
\end{example}

\section{The Interpretation with Reflection Principle}
\label{intrefl}

Let $\LLL$ be the language of first-order predicate logic of signature $\Sigma$ as defined in Section \ref{indefappsec}.

\subsection{Defining the Interpretation.}
\label{intinreflp}

In the previous two sections we introduced the requirements in order to define the interpretation $\MMM_\III \validb{\yA}{\vecb{a}:C}$ of a formula $\yA$ in an indefinitely extensible structure $\MMM_\III$. We often omit the structure, writing $\validb{\yA}{\vecb{a}:C}$ for $\MMM_\III \validb{\yA}{\vecb{a}:C}$. The interpretation uses the following concepts:
\begin{enumerate}
\item A signature $\Sigma^\III$ over $\Sigma$ and a $\Sigma^\III$-structure $\MMM_\III$, see Section \ref{unimultsec}.
\item A $\ll$-relation $\ll$, see Section \ref{relinf}.
\item A state declaration $\tdm{C}{\yA}$, see Section \ref{typedecltermform}.
\item A variable assignment $\vecb{a} = (a_0,\dots, a_{n-1}) \in \MMM_C$.
\end{enumerate}

The variable assignment $\vecb{a}$ is thus relative to the context $C$, given by the state declaration. The structure $\MMM_\III$ together with the interpretation yields the $\Sigma^\III$-model $(\MMM_\III,\validbn)$.

\begin{definition}
\label{termdef}
Given a $\Sigma^\III$-structure $\MMM_\III$. The interpretation of variables and formulas are defined recursively on the state declaration. The interpretation $\val{\ }$ is a function with two arguments, first, a variable in a state context $\tdt{\yx{k}}{C}{j}$, i.e.~the triple $(C,\yx{k},j)$, and second, a variable assignment $\vecb{a} = (a_0,\dots, a_{n-1}) \in \MMM_C$. Its value is in $\MMM_j$. The interpretation of a formula $\validbn$ is a relation between a formula in a state context $\tdm{C}{\yA}$ and a variable assignment $\vecb{a} = (a_0,\dots, a_{n-1}) \in \MMM_C$.
\begin{eqnarray*}
\valp{\yx{k}}{j}_{\vecb{a} : C} &:=& a_k \ \text{ for } C = (i_0, \dots,i_{n-1}) \text{ and } \tdt{\yx{k}}{C}{j}.\\
\validb{\y{R} \yt{0} \dots \yt{m-1}}{\vecb{a}:C} &\iffdef& R_{(j_0,\dots,j_{m-1})}(\valp{\yt{0}}{j_0}_{\vecb{a} : C}, \dots , \valp{\yt{m-1}}{j_{m-1}}_{\vecb{a} : C}) \\
&&\text{for some } \tdr{\y{R}}{(j_0,\dots, j_{m-1})}.\\
\validb{\yA \to \yB}{\vecb{a}:C} &\iffdef& \validb{\yA}{\vecb{a}:C} \text{ implies } \validb{\yB}{\vecb{a}:C},
\end{eqnarray*}
and accordingly for the connectives $\land$ and $\lor$. $\validb{\bot}{\vecb{a}:C}$ holds for no $\vecb{a} \in \MMM_C$. The interpretation of the quantifiers is:
\begin{eqnarray*}
\validb{\forall \y{x} \yA}{\vecb{a}:C} &\iffdef& \validb{\yA}{\vecb{a}b:Ci} \text{ holds for all elements}  \\
&& b \in \MMM_i \ \text{ for some } i \in \III \text{ (with } C \ll i).\\
\validb{\exists \y{x} \yA}{\vecb{a}:C} &\iffdef& \validb{\yA}{\vecb{a}b:Ci} \text{ holds for an element} \\
&& b \in \MMM_i \ \text{ for all } i \in \III \text{ (with } C \ll i).
\end{eqnarray*}

Moreover, $\val{\yA}$ denotes the family $(\val{\yA}_C)_{C \in \HHH}$ with $\HHH = \{C \in \III^n \mid \tdm{C}{\yA}\}$, such that $\val{\yA}_C \subseteq \MMM_C$ and $\vecb{a} \in \val{\yA}_C \, \iffdef \validb{\yA}{\vecb{a}:C}$. A $\Sigma^\III$-model $(\MMM_\III,\validbn)$ is called \emph{indefinitely large} iff the signature $\Sigma^\III$ and the $\ll$-relation are indefinitely large.
\end{definition}

The condition $C \ll i$ in brackets (used for the quantifiers) is already a consequence of the declaration $\tdm{C}{\yA}$ and can be omitted. The interpretation $\val{\neg \yA}$ is the family $(\MMM_C \setminus \val{\yA}_C)_{\tdm{C}{\yA}}$, so negation is interpreted as a family of relative complements. 

It is important to note that Definition \ref{termdef} is at this stage too general as to be useful. The interpretation of a formula $\yA$ depends by definition on the state declaration of $\yA$ and could be different for different declarations. It is vital to have suitable restrictions on the relation $\ll$, call \emph{adequate} (see Definition \ref{adeqdef}). The next example is one of many which shows the provisional nature of the definition so far.

\begin{example}
\label{pathoex}
Consider the model $(\nat_\nat,\validbn)$ and let $i \gg (i_0,\dots,i_{n-1})$ hold iff $i \geq i_0$ and $\dots$and $i \geq i_{n-1}$. Then we have $\nat_\nat \validb{(\forall \yx{1} \, \yx{1} \leq \yx{0})}{0:1}$, since there is an index $i \gg 1$ such that $b \leq 0$ for all $b \in \nat_i$ (take $i=1$). But $\nat_\nat \not\validb{(\forall \yx{1} \, \yx{1} \leq \yx{0})}{0:2}$ since $b \leq 0$ is not satisfied for all elements in $\nat_i$ if $i \geq 2$ (take $1$ for $\yx{1}$).
\end{example}

In contrast to the universal quantifier, the interpretation of the existential quantifier uses the locution ``for all $i \in \III$'', so $i$ does not range over a fixed finite collection\footnote{We use this formulation mainly to be in line with the classical relation between existential and universal quantifiers.}. In Section \ref{adeqint} we will see that it possible to replace ``for all $i \in \III$'' by a single index $i \in \III$. It is easy to confirm that the interpretation $\validbn$ becomes the usual Tarskian interpretation $\valid$ if the system consists of one single set $\MMM_i$. Then $i^n$ must be a $\ll$-context for all $n \in \nat$ and the structure must be finite then.

\begin{remark}
An open formula $\yA(\yx{0},\dots,\yx{n-1})$ is implicitly seen as universally quantified, that is, it corresponds to the sentence $\forall \yx{0} \dots \yx{n-1} \yA$. For the common interpretation in a Tarskian model $(\MMM,\valid)$, this leads to $\valid \forall \yx{0} \dots \yx{n-1} \yA \iff \validv{\yA}{\vecb{a}}$ for all variable assignments $\vecb{a} \in \MMM^n$. For the interpretation $\validbn$ however we have $\validbn \forall \yx{0} \dots \yx{n-1} \yA \iff \validb{\yA}{\vecb{a} : C}$ for all variable assignments $\vecb{a} \in \MMM_C$ of \emph{some} $\ll$-context $C$.
\end{remark}

\begin{lemma}
\label{indeflmodlem}
Each definable relation in an indefinitely large model $(\MMM_\III,\validbn)$ is indefinitely large.
\end{lemma}

\begin{proof}
Each definable relation is of the form $\val{\yA}$. The set of state declarations for $\yA$ is indefinitely large by Lemma \ref{apxtypdecllem2}, hence $\val{\yA}$ is indefinitely large, too, since the index set of $\val{\yA}$ comprises the contexts $C$ with $\tdm{C}{\yA}$ by definition.
\end{proof}

\subsection{How to Find Sufficiently Large Indices.}
\label{sufflarge}

Our aim is to show that propositions $\yA$ interpreted by $\validbn$ have the same truth value as in the usual Tarskian model and are thus independent of the chosen state declaration. Additionally, the value $\val{\yA}$ will then be a (compatible) relation in the sense of Definition \ref{dircomplemma}. The technique which we apply for this purpose has been used in various ways, most notably in the L\"owenheim-Skolem theorem. Basically, we must guarantee that an index $i \gg C$ is as large as to embrace all witnesses of valid existential quantified formulas in scope.

We may describe the notion $C \ll i$ also as the use of ``means''. If $C$ defines the current stage, then the elements of consideration are inside $\MMM_C$, in particular all assignments $\vecb{a}$ are within $C$. If we claim $\forall \yx{n}\,\yA(\yx{0},\dots,\yx{n})$, then we only have to make sure that we consider all elements $b$, replacing $\yx{n}$, that are reachable by the means applied to the current stage $\MMM_C$, with $a_0 \in \val{i_0}, \dots, a_{n-1} \in \val{i_{n-1}}$ replacing $\yx{0},\dots,\yx{n-1}$ for $C = (i_0,\dots,i_{n-1})$. These means are given by a set of relations such that at each stage there is only a finite set $\SSS$ of them. All relations in $\SSS$ are moreover definable from a finite set $\TTT \subseteq \LLL$ of formulas\footnote{If $\LLL$ contains function symbols, then we have to add terms to $\LLL$ and functions to $\SSS$.}. So the construction only needs the definable relations. 

This yields the following chain --- we add the subscript $\TTT$ to $\SSS$ and $\ll$ in order to indicate the dependency from the set $\TTT$:
\begin{align*}
\text{Set of formulas } \TTT \ &\leadsto \ \text{Set of relations } \SSS_\TTT := \{\val{\yB} \mid \exists \yx{} \yB \in \hat\TTT\} \\ 
&\leadsto \ \text{Relation } \ll_\TTT \\
&\leadsto \ \text{Interpretation } \validbn.
\end{align*}

To get $\hat\TTT$ from $\TTT$, first replace each occurrence of $\forall \y{x}$ in $\TTT$ by $\neg \exists \y{x} \neg$. Then add all subformulas resulting in the set $\hat\TTT$, which is finite whenever $\TTT$ is. We will introduce the other notions soon.

\subsection{Avoiding a Circularity.}
\label{avcirc}

The step from a set of formulas $\TTT$ to the set of relations $\SSS_\TTT = \{\val{\yB} \mid \exists \yx{} \yB \in \hat\TTT\}$ already requires an interpretation. The natural choice would be to take the interpretation $\validbn$. But then we already need an adequate relation $\ll$ before we define the properties that we require from $\ll$. In order to avoid this circularity, define an auxiliary interpretation $\validun$, which does not use relation $\ll$ and which yields the same truth values as $\validbn$ for formulas in $\TTT$.

Let the state declaration be that from Definition \ref{typedecldef} but without condition $i \gg C$ in case of the quantifiers, i.e., we take the trivial $\ll$-relation $\III^*$. It is easy to confirm that for an indefinitely large signature $\Sigma^\III$ the assignment $\tdm{C}{\yA}$ holds for any context $C$ of length $n$ and any $n$-ary formula $\yA$ (the proof is similar as for Lemma \ref{apxtypdecllem2}, but simpler). In other words, each formula $\yA$ is approximable in $\Sigma^\III$ with $\III^*$. The interpretation $\validun$ is as $\validbn$, but with the quantifiers adopted as follows:
\begin{eqnarray*}
\validu{\forall \y{x} \yA}{\vecb{a}:C} &\iffdef& \validu{\yA}{\vecb{a}b:Ci}\text{ for all } b \in \MMM_i \text{ for all } i \in \III.\\
\validu{\exists \y{x} \yA}{\vecb{a}:C} &\iffdef& \validu{\yA}{\vecb{a}b:Ci}\text{ for some } b \in \MMM_i \text{ for some } i \in \III.
\end{eqnarray*}

So interpretation $\validun$ simply simulates on an increasing carrier set the Tarskian interpretation. The relation generated by the formulas $\yA$ is denoted as $\val{\yA}^m$, that is $\vecb{a} \in \val{\yA}^m_C \, \iffdef \validu{\yA}{\vecb{a}:C}$. The proof of the next lemma is straightforward.

\begin{lemma}
\label{unimultthm2}
Given a $\Sigma^\III$-structure $\MMM_\III$ with an indefinitely large signature $\Sigma^\III$. Then $(\MMM_\III,\validun)$ is an indefinitely large model. For each formula $\yA$ in $\LLL$ we have $\MMM_\III \validu{\yA}{\vecb{a}:C} \iff \bigcup \MMM_\III \validv{\yA}{\vecb{a}}$ for all $\vecb{a} \in \MMM_C$.
\end{lemma}

The index set of an $n$-ary relation $\val{\yA}^m$ is not only in $\fil_n$, but it is $\III^n$ (provided the signature $\Sigma^\III$ is indefinitely large). Moreover, $\validu{\yA}{\vecb{a}:C} \iff \validu{\yA}{\vecb{a}:C'}$ holds for all assignments $\vecb{a} \in \MMM_C \cap \MMM_{C'}$.

\subsection{Witnesses of Existential Quantified Formulas.}
\label{defmeans}

Given an $n+1$-ary relation $R_\HHH$ and an element $\vecb{a} \in \MMM_C$ for a context $C$ of length $n$. The following sets play a key role:
\begin{equation}
\label{irsetdef}
\iset{R} := 
\begin{cases}
\{i \in \III \mid \exists b \in \MMM_i \ R_{Ci}(\vecb{a}b)\} & \text{if such $i$ exists,}\\
\III & \text{otherwise.}
\end{cases}
\end{equation}

This set is in any case non-empty. In a next step we define a $\ll$-relation $\ll_R$ from this set. First define $(\ll_R)_{n+1}$ by
\begin{equation}
\label{muupper}
C \ll_R i \ \iffdef i \in \bigcap_{\vecb{a} \in \MMM_C} \iset{R},
\end{equation}
and let $\ll_R$ be the $\ll$-relation generated by $(\ll_R)_{n+1}$ (see Definition \ref{gendef}). Note that $\ll_R$ is indeed a $\ll$-relation by Lemma \ref{extlem}, i.e., $Ci \in (\ll_R)_{m+1}$ implies $C \in (\ll_R)_m$ for a context $C$ of length $m$.

\begin{lemma}
\label{ggsetupwclr}
If $R_\HHH$ is an indefinitely large $n+1$-ary relation, then $\ll_R$ is an indefinitely large $\ll$-relation.
\end{lemma}

\begin{proof}
By Lemma \ref{extlem} it suffices to show $(\ll_R)_{n+1} \in \fil_{n+1}$, that is 
\begin{equation}
\label{preq}
\{C \in \III^n \mid \bigcap_{\vecb{a} \in \MMM_C} \iset{R} \in \fil_1 \} \in \fil_n.
\end{equation}

By assumption we have $\{C \in \III^n \mid C^\HHH \in \fil_1\} \in \fil_n$ since $\HHH \in \fil_{n+1}$. Furthermore, $\up i \cap C^\HHH \subseteq \iset{R}$ holds for any $\vecb{a} \in \MMM_C$ and $i \in \iset{R}$: In the first case of definition (\ref{irsetdef}), $R_{Ci}(\vecb{a}b)$ implies $R_{Ci'}(\vecb{a}b)$ for all $i' \geq i$ with $Ci' \in \HHH$ by compatibility of $R_\HHH$; for the ``otherwise case'' this holds trivially.

Therefore, if $C^\HHH \in \fil_1$ then $\iset{R}$ is also in $\fil_1$ and we conclude that $\{C \in \III^n \mid \iset{R} \in \fil_1\} \in \fil_n$. This holds for all finitely many $\vecb{a} \in \MMM_C$, hence Property (\ref{preq}) follows from the fact that $\bigcap_{\vecb{a} \in \MMM_C} \iset{R} \in \fil_1$ holds iff $\iset{R} \in \fil_1$ for all $\vecb{a} \in \MMM_C$.
\end{proof}

To a set $\SSS$ of relations on a system $\MMM_\III$, define $\ll_\SSS$ by
\begin{equation*}
C \ll_\SSS i \ \iffdef \ C \ll_R i \text{ for all } R \in \SSS.
\end{equation*}

It follows immediately from Lemmata \ref{ggsetupwclr} and \ref{indefcor}:

\begin{corollary}
\label{ggsetupwcl}
Given a finite set $\SSS$ of indefinitely large relations $R_\HHH$, then $\ll_\SSS$ is an indefinitely large $\ll$-relation.
\end{corollary}

With the interpretation\footnote{If one is willing to accept the circularity, then it is possible to use interpretation $\validbn$ instead of $\validun$.} $\validun$ from Section \ref{avcirc} define 
\[
\iset{\exists \yx{} \yB} := \iset{\val{\yB}^m} \quad \text{and} \quad \SSS_\TTT := \{\val{\yB}^m \mid \exists \yx{} \yB \in \hat\TTT\}.
\]

Remind that $\hat\TTT$ results from $\TTT$ by replacing each occurrence of $\forall \y{x}$ in $\TTT$ by $\neg \exists \y{x} \neg$ and then adding all subformulas. Let the relation $\ll_{\SSS_\TTT}$ be abbreviated by $\ll_\TTT$.

\begin{example}
Let $\TTT$ consist of the single formula $\exists \yx{1} \y{R} \yx{0} \yx{1}$ and consider the structure $\nat_\nat$. The interpretation of $\y{R}$ shall be ``$\yx{0}+\yx{1}$ is a perfect number'', the first of them are 6 and 28. The relevant relation in $\SSS_\TTT$ is $R_\HHH$ with $\HHH = \nat \times \nat$ and $R_{(i_0,i_1)} := \{(a_0,a_1) \in \nat_{i_0} \times \nat_{i_1} \mid P_{i_0 + i_1 -1}(a_0 + a_1)\}$, where $P_i(a)$ holds iff $a < i$ is a perfect number. Consider the context $C = 8$. We have $\isetp{R}{a_0 : i_0} = \{i_1 \in \nat \mid \exists a_1 \in \nat_{i_1} \, P_{i_0 + i_1 -1}(a_0 + a_1)\}$, thus
\[
\bigcap_{a_0 \in \nat_8 = \{0,\dots,7\}} \isetp{R}{a_0:8} = \ \up 7 \ \cap \up 6 \ \cap \dots \cap \up 1 \ \cap \up 22 = \ \up 22,
\]
since the sets $\up 7$ to $\up 1$ stem from assigning $0$ up to $6$ to $a_0$, while $\up 22$ is needed for $a_0 = 7$ (note that $P_{29}(7 + 21)$ holds and $21 \in \nat_{22}$). Hence we have $8 \ll_\TTT 22$ and $(8,22)$ is thus a $\ll$-context.
\end{example}

\begin{remark}
\label{axrem}
If we consider a set of valid sentences $\TTT$, as for instance axioms, then only the restricted occurrences are relevant, i.e., the positive occurrences of $\exists \yx{}$ and the negative occurrences of $\forall \yx{}$, since for other occurrences no witness (or counterexample) exists. So instead of replacing each occurrence of a universal quantifier to yield $\hat\TTT$, it is enough in that situation to only translate negative occurrences of $\forall \yx{}$.
\end{remark}

\subsection{Adequacy.}
\label{adeqint}

We are now ready to define the necessary restriction on $\ll$-relations in order to get a correct interpretation. This section contains the central result, Proposition \ref{coincidethm}, from which the final theorems follow easily.

\begin{definition}
\label{adeqdef}
A $\ll$-relation $\ll$ is \emph{adequate} for a set of formulas $\TTT$ iff it is a subset of $\ll_\TTT$. An interpretation $\validbn$, and similarly a model $(\MMM_\III,\validbn)$, is called adequate for $\TTT$ if the underlying relation $\ll$ is adequate for $\TTT$.
\end{definition}

A simple consequence of this definition is that if a $\ll$-relation on $\III$ is adequate for $\TTT$, then its restriction to a subset $\JJJ \subseteq \III$ is also adequate for $\TTT$.

\begin{lemma}
\label{adexlem}
Given an indefinitely large signature $\Sigma^\III$ and a finite set $\TTT$ of formulas. Then there is an indefinitely large $\ll$-relation that is adequate for $\TTT$.
\end{lemma}

\begin{proof}
Take the relation $\ll_\TTT$, which is adequate for $\TTT$, and use Corollary \ref{ggsetupwcl} to show that $\ll_\TTT$ is indefinitely large. Firstly, $\SSS_\TTT$ is a finite set of relations since $\TTT$ is finite. Secondly, all relations in $\SSS_\TTT$ are indefinitely large (as a consequence of Lemmata \ref{indeflmodlem} and \ref{unimultthm2}).
\end{proof}

Note that this is different to a situation in which a relation $\ll$ exists which is adequate for the set $\LLL$ of \emph{all} expressions. For instance, for a model $(\nat_\nat,\validbn)$ of arithmetic there is no such $\ll$-context except $\nil$. This is due to the fact that $\LLL$ contains the valid formulas $\exists \yx{0} \, \yx{0} \geq \y{n}$ for each $n \in \nat$, and $\bigcap_{n \in \nat} \III_{\exists \yx{0} \, \yx{0} \geq \y{n}} = \emptyset$. 

The next proposition states that within $\TTT$ the interpretations $\validbn$ and $\validun$ are the same and they coincide with the usual validity in Tarskian semantics. Let $\bigcup \MMM_\III$ be the structure with underlying set $\bigcup_{i \in \III} \MMM_i$, defined in Section \ref{unimultsec}.

\begin{proposition}
\label{coincidethm}
Given a model $(\MMM_\III,\validbn)$, adequate for a set $\TTT$ of formulas. Then for each $\yA \in \TTT$, each $\ll$-context $C$ with $\tdm{C}{\yA}$ and each variable assignment $\vecb{a} \in \MMM_C$ we have
\[
\MMM_\III \validb{\yA}{\vecb{a} : C} \ \iff\ \MMM_\III \validu{\yA}{\vecb{a} : C} \iff \bigcup \MMM_\III \validv{\yA}{\vecb{a}}.
\]
\end{proposition}

\begin{proof}
Consider in the beginning the situation that $\TTT$ is closed under sub-expressions and does not contain a universal quantifier. The equivalence is shown by induction on $\tdm{C}{\yA}$ and the only interesting case is the existential quantifier $\yA = \exists \y{x} \yB$. The implication ``$\validbn \ \imp\ \validun$'' follows from the fact that for all $\ll$-contexts $C$ there is an index $i \gg C$, since the declaration $\tdm{C}{\exists \y{x} \yB}$ presupposes $\tdm{Ci}{\yB}$ with some $i \gg C$. Next we show the implication ``$\validun \ \imp\ \validbn$'', i.e.,
\begin{multline*}
\text{there is } j \in \III \text{ such that } \exists b \in \MMM_j \ \validu{\yB}{\vecb{a}b : Cj} \\
\text{implies }\exists b' \in \MMM_{i'} \ \validb{\yB}{\vecb{a}b' : Ci'} \text{ for all } i' \gg C.
\end{multline*}

Given $\vecb{a} \in \MMM_C$. Formula $\exists \y{x} \yB$ is an element of $\TTT = \hat\TTT$, hence $R := \val{\yB}^m$ is an element of $\SSS_\TTT$. We must prove $\exists b' \in \MMM_{i'} \ \validb{\yB}{\vecb{a}b' : Ci'}$ for all $i' \gg C$. By induction hypothesis this is the same as $\exists b' \in \MMM_{i'} \, R_{Ci'}(\vecb{a}b')$ for all $i' \gg C$. 

Every index $i'$ with $C \ll i'$ satisfies $C \ll_\TTT i'$ since $\validbn$ is adequate for $\TTT$, hence $C \ll_R i'$. Consequently, $i' \in \iset{R}$. By assumption there is an index $j \in \III$ such that $\exists b \in \MMM_j \, R_{Cj}(\vecb{a}b)$. Therefore, by definition of set $\iset{R}$, there is an element $b' \in \MMM_{i'}$ with $R_{Ci'}(\vecb{a}b')$. 

So we are finished for the case that $\TTT$ is closed under sub-expressions and does not contain a universal quantifier. The general situation is reduced to the just proven one by noticing that $\ll_\TTT$ is the same as $\ll_{\hat\TTT}$ (or by using the fact that $\forall \y{x}$ can be reduced to $\neg \exists \y{x} \neg$ for both interpretations). The second equivalence has been stated in Lemma \ref{unimultthm2}.
\end{proof}

This proposition has useful consequences. First of all it shows that an interpretation $\validbn$, which is adequate for $\TTT$, does not diverge from the common interpretation in a Tarskian model, as long as formulas are taken from the set $\TTT$. For formulas outside $\TTT$ there will surely be differences. From a potentialist's point of view however, $\TTT$ is at most potential infinite, and the substitute for talking about the whole set $\TTT$ is to specify a sufficiently large finite set of expressions, so there is no ``outside''.

Secondly, the interpretation $\validbn$ is independent of the chosen state declaration and the way it has been derived. This follows from the equivalence $\MMM_\III \validb{\yA}{\vecb{a} : C} \iff \bigcup \MMM_\III \validv{\yA}{\vecb{a}}$. Similarly we have:

\begin{corollary}
\label{compcor}
Given a model $(\MMM_\III,\validbn)$, adequate for a set $\TTT$ of formulas, then $\val{\yA}$ is a (compatible) relation for all approximable formulas $\yA \in \TTT$.
\end{corollary}

\begin{proof}
Compatibility of $\val{\yA}$ is a consequence of the following equivalences $\MMM_\III \validb{\yA}{\vecb{a} : C} \iff \bigcup \MMM_\III \validv{\yA}{\vecb{a}} \iff \MMM_\III \validb{\yA}{\vecb{a} : C'}$ for $\vecb{a} \in \MMM_C \cap \MMM_{C'}$. Since $\yA$ is approximable, i.e., $\tdm{C}{\yA}$ for some context $C$, there is at least one instance $\val{\yA}_C$, hence $\val{\yA}$ is a relation.
\end{proof}

As a further consequence from Proposition \ref{coincidethm} we have an alternative definition for $\validb{(\exists \y{x} \yA)}{\vecb{a}:C}$, that is, $\validb{\yA}{\vecb{a}b:Ci}$ holds for an element $b \in \MMM_i$ for \emph{some} $i \in \III$. This follows immediately from Proposition \ref{coincidethm} and the interpretations of $\exists \y{x} \yA$ with respect to $\validun$ and $\validbn$ resp. Basically this states the irrelevance of the outer quantifier, be it $\exists i \gg C$ or $\forall i \gg C$ or a fixed $i \gg C$. This is indeed an essential property of an index being ``sufficiently large''. The next corollary expresses this irrelevance of the index $i \gg C$:

\begin{corollary}
\label{allexgen}
Given a set $\TTT$ of formulas with $\yA \in \TTT$, a model $(\MMM_\III,\validbn)$ adequate for $\TTT$, and a $\ll$-context $C$ with $\tdm{C}{\yA}$. Then for any $i \gg C$ and any variable assignment $\vecb{a} \in \MMM_C$ we have
\begin{align*}
\validb{\forall \y{x} \yA}{\vecb{a}:C} &\iff \validb{\yA}{\vecb{a}b:Ci} \text{ holds for all elements } b \in \MMM_i.\\
\validb{\exists \y{x} \yA}{\vecb{a}:C} &\iff \validb{\yA}{\vecb{a}b:Ci} \text{ holds for an element } b \in \MMM_i.
\end{align*}
\end{corollary}

\subsection{Soundness and Completeness.}
\label{soundcomp2}

In this section we show that the constructions of $\Sigma$-models out of $\Sigma^\III$-models and vice versa preserve validity. As a consequence, the usual deduction rules of classical first-order predicate logic are sound and complete with respect to the collection $\allapp$ of all (indefinitely extensible) $\Sigma^\III$-models. Thereby we not only have to vary\footnote{Consider the sentence $\yA := \forall \y{x} \, \y{P} \y{x}$ and the models $(\nat,\valid)$ and $(\nat_\nat,\validbn)$. The formula $\yA$ is valid in $\nat$ iff $P(b)$ holds for all $b \in \nat$, whereas its validity in $\nat_\nat$ depends on the $\ll$-relation. For a fixed $\ll$-relation, formula $\yA$ is valid iff $P(b)$ holds for all $b \in \nat_{i_0}$, for some $i_0 \gg \nil$, so validity in both kinds of model differ.} on the set of models, but also on the index set $\III$ and relation $\ll$.

Let $\TTT$ be a set of sentences (of a language $\LLL$ of signature $\Sigma^\III$), which is no restriction since we may consider the universal closure of an open formula. We write $\allapp \validbn \TTT$ to mean that $\MMM_\III \validbn \yA$ holds for all formulas $\yA \in \TTT$ and all $\Sigma^\III$-models $(\MMM_\III,\validbn) \in \allapp$, adequate for $\TTT$. Note that in $\allapp \validbn \TTT$, the suffix $\ll$ is ``generic'', not referring to a specific relation.

\begin{proposition}
\label{soundcompllem}
Let $\TTT$ be a set of sentences and $\yA \in \TTT$.
\begin{enumerate}
\item For each $\Sigma^\III$-model $(\MMM_\III,\validbn)$, adequate for $\TTT$, there is a Tarskian $\Sigma$-model $(\MMM,\valid)$ such that $\MMM_\III \validbn \yA \iff \MMM \valid \yA$.

\item If $\TTT$ is finite, then for each Tarskian $\Sigma$-model $(\MMM,\valid)$ there is an indefinitely large $\Sigma^\III$-model $(\MMM_\III,\validbn)$ which is adequate for $\TTT$ and which satisfies $\MMM_\III \validbn \yA \iff \MMM \valid \yA$.
\end{enumerate}
\end{proposition}
 
\begin{proof}
In order to show the first claim, take the union $\MMM = \bigcup \MMM_\III$ and apply Proposition \ref{coincidethm}. For the second claim consider for a given structure $\MMM$ the structure $\MMM_\III$ with $\III \simeq \potfin{\MMM}$, defined in Section \ref{unimultsec}, which satisfies $\MMM = \bigcup \MMM_\III$. Each relation has an indefinitely large index set $\HHH = \III^n$, and $\Sigma^\III$ is indefinitely large, both by Lemma \ref{multlem}. Because $\Sigma^\III$ is indefinitely large we may apply Lemma \ref{adexlem} to get an indefinitely large relation $\ll$, adequate for $\TTT$. Finally, apply Proposition \ref{coincidethm} again.
\end{proof}

Let $\stdmod$ be the collection of all usual Tarskian models for a given signature $\Sigma$ and $\stdmod \valid \TTT$ denote $\MMM \valid \yA$ for all $(\MMM,\valid) \in \stdmod$ and $\yA \in \TTT$. The following corollary follows immediately from the proposition above.

\begin{corollary}
\label{soundcomplthm1}
Given a finite set $\TTT$ of sentences, then
\[
\allapp \validbn \TTT \iff \stdmod \valid \TTT.
\] 
\end{corollary}

This is obviously the basis of a soundness and completeness result by applying the corresponding theorems for Tarskian models\footnote{This is an indirect argument to show completeness. Moreover, we have shown completeness based on actual infinite Tarskian models. A potentialist might not accept this (see however Section \ref{reflsec}). Though it is possible to prove the result directly. A typical construction of a model from a consistent set of formulas, e.g., that of Henkin, already gives a model with an increasing family as underlying carrier set. This is straightforward for theories without equality, for a theory with equality, the construction of the term model however requires a non-injective embeddings between the sets $\MMM_i$ --- see Section \ref{simpadopt}.}. From a potentialist's viewpoint, either $\TTT$ is a fixed finite set and we apply Corollary \ref{soundcomplthm1} to it. Or, if the set $\TTT$ is potential infinite, we use some (indefinitely large) finite subset of it and again apply Corollary \ref{soundcomplthm1}.

\begin{corollary}
\label{soundcompl2}
From a potentialist's viewpoint, the interpretation $\validbn$ is sound and complete with respect to the collection $\allapp$ of all indefinitely extensible models and a common deductive system of classical first-order predicate logic.
\end{corollary}

\subsection{A Finite Submodel.}
\label{smallsub}

Let a set of formulas $\TTT$ be given as well as a $\Sigma^\III$-model $(\MMM_\III,\validbn)$, adequate for $\TTT$. There is a submodel of $\MMM_\III$ that suffices to interpret the formulas in $\TTT$ correctly. This submodel is finite, whenever $\TTT$ is.  Let us call a structure $\MMM_\JJJ$ of signature $\Sigma^\JJJ$ a \emph{substructure} of a $\Sigma^\III$-structure $\MMM_\III$ iff $\JJJ \subseteq \III$ and $\Sigma^\JJJ$ consists of those $\tdr{\y{R}}{C}$ in $\Sigma^\III$, for which $C \in \JJJ^n$ holds. That is, the interpretation of $\y{R} \in \Sigma$ as a relation $R_\HHH$ in $\MMM_\JJJ$ is a subset of instances from the interpretation in $\MMM_\III$.

\begin{definition}
Given a set of formulas $\TTT$. A model $(\MMM_\JJJ,\validbn)$ is a \emph{$\TTT$-submodel} of $(\MMM_\III,\validbn)$ iff $\MMM_\JJJ$ is a substructure of $\MMM_\III$ and 
\begin{equation}
\label{elemsubstr}
\MMM_{\JJJ} \validb{\yA}{\vecb{a} : C} \iff \MMM_{\III} \validb{\yA}{\vecb{a} : C}
\end{equation}
holds for all $\yA(\yx{0}, \dots, \yx{n-1}) \in \TTT$, all $\ll$-contexts $C \in \JJJ^n$ with $\tdm{C}{\yA}$, and all assignments $\vecb{a} \in \MMM_C$. 
\end{definition}

If $(\MMM_\III,\validbn)$ is a $\Sigma^\III$-model and $\JJJ$ is a possible restriction for $\TTT$ (see Definition \ref{sufflargedef}), then each formula approximable in $\Sigma^\III$ is also approximable in $\Sigma^\JJJ$ and thus has an interpretation $\validbn$ in the structure $\MMM_\JJJ$. The $\ll$-relation used for the interpretation in $\MMM_\JJJ$ is the restriction of that in $\MMM_\III$ to $\JJJ$. We will first show that a substructure $\MMM_{\JJJ}$ of $\MMM_\III$ is automatically a $\TTT$-submodel in that case. 

\begin{lemma}
\label{finsublem}
Given a set of formulas $\TTT$, approximable in $\Sigma^\III$ with $\ll$. Let $\ll$ be adequate for $\TTT$ and let $\JJJ \subseteq \III$ be a possible restriction for $\TTT$. Then $(\MMM_\JJJ,\validbn)$ is a $\TTT$-submodel of $(\MMM_\III,\validbn)$, which is adequate for $\TTT$.
\end{lemma}

\begin{proof}
We show (\ref{elemsubstr}) by induction on the state declaration $\tdm{C}{\yA}$. Since the satisfaction relation is independent of the state declaration and the way it has been derived (see Section \ref{adeqint}), we may assume that the same declaration $\tdm{C}{\yA}$ (with the same derivation) as for $\MMM_\JJJ$ has been used for $\MMM_\III$ as well. 

This immediately yields the equivalence for atomic formulas $\y{R} \yt{0} \dots \yt{m-1}$. For $\yA \to \yB$, $\yA \land \yB$, and $\yA \lor \yB$ the claim follows straightforwardly from the induction hypothesis, so consider $\exists \y{x} \yB$. By induction hypothesis we have for all $\vecb{a}b:Ci$, with the index $i \in \JJJ$, $i \gg C$, used in the derivation $\tdm{Ci}{\yB}$:
\[
\MMM_{\JJJ} \validb{\yB}{\vecb{a}b:Ci} \iff \MMM_\III \validb{\yB}{\vecb{a}b:Ci}.
\]

We can use Corollary \ref{allexgen} on both sides of the equivalence to get the equivalence $\MMM_{\JJJ} \validb{(\exists \y{x} \yB)}{\vecb{a}:C} \iff \MMM_{\III} \validb{(\exists \y{x} \yB)}{\vecb{a}:C}$. A similar consideration applies to universal quantified formulas.
\end{proof}

\begin{theorem}
\label{finiterestthm}
Given a finite set $\TTT$ of formulas as well as a finite set of indices $\III_0 \subseteq \III$. Let the $\Sigma^\III$-model $(\MMM_\III, \validbn)$ be indefinitely large and adequate for $\TTT$. Then there is a finite set $\JJJ$ with $\III_0 \subseteq \JJJ \subseteq \III$ such that $(\MMM_\JJJ,\validbn)$ is a $\TTT$-submodel of $(\MMM_\III,\validbn)$, which is adequate for $\TTT$.
\end{theorem}

\begin{proof}
$\TTT$ is approximable in $\III$ with $\ll$ by Lemma \ref{apxtypdecllem2}. By Lemma \ref{finiterestlem} there is a possible restriction $\JJJ \subseteq \III$ for $\TTT$ with $\III_0 \subseteq \JJJ$. So the result follows from Lemma \ref{finsublem}.
\end{proof}

Consider the situation of Example \ref{simpex} with a finite set of formulas $\TTT$. The finite set of indices $\{i_0, \dots i_{n-1}\}$, mentioned at the end of the example, suffice to define the $\TTT$-submodel. More precisely, the index set is $\JJJ = \{i_0, \dots i_{n-1},j\}$ with upper bound $j$, if necessary. The interpretation of the sentence $\yA := \forall \yx{0} \, \exists \yx{1} \, \forall \yx{2} \, \y{R}(\yx{0},\yx{1},\yx{2})$, for instance, is
\begin{equation*}
\validbn \yA \iff \forall a_0 \in \MMM_{i_0} \, \exists a_1 \in \MMM_{i_1} \, \forall a_2\in \MMM_{i_2} \, R_{(i_0,i_1,i_2)}(a_0,a_1,a_2),
\end{equation*}
which is the same in $\MMM_\III$ as in $\MMM_\JJJ$.

\subsection{Extended Example: Set Theory.}
\label{extexsec}

Shaughan Lavine considered in \cite{lavine2009understanding} ZFC set theory and called the finitistic translation of the axiom of infinity (by adding bounds to the variables) ``Axiom of a Zillion''. In our approach the axiom of a zillion is the usual axiom of infinity, but interpreted in an increasing model. The signature $\Sigma$ of the language of set theory has two binary relation symbols $\eps$ (membership) and $=$ (equality).

Consider a set universe\footnote{We may either assume that $\VVV$ is the class of all sets, or it is already an increasing structure, see the remarks in Section \ref{reflsec}.} $\VVV$ with index set $\III = \potfin{\VVV}$ and $\VVV_i = i$. The signature $\Sigma^\III$ over $\Sigma$ consists of $\eps : (i_0,i_1)$ and $= \, : (i_0,i_1)$ for all $i_0,i_1 \in \III$. We use common abbreviations, e.g., $0 := \emptyset$, $1 := \{\emptyset\}$, $2 := \{\emptyset,\{\emptyset\}\}$, $\yA \subseteq \yB := \forall \yx{} (\yx{} \eps \yA \to \yx{} \eps \yB)$,
\begin{align*}
\yx{} \eps \yA \setminus \yB &:= \yx{} \eps \yA \land \neg \yx{} \eps \yB, \\
\yx{} \eps \yA \triangle \yB &:= \yx{} \eps \yA \setminus \yB \lor \yx{} \eps \yB \setminus \yA, \\
\yx{1} = \y{suc} \, \yx{0} &:= \forall \yx{2} ( \yx{2} \eps \yx{1} \yequi \yx{2} = \yx{0} \lor \yx{2} \eps \yx{0}).
\end{align*}

Let the ZFC axioms be given in an enumerated form. As mentioned in Remark \ref{axrem}, it suffices to avoid only the negative occurrences of universal quantifiers in order to yield set $\hat{\TTT}$. Let the first 5 axioms --- already formulated without negative occurrences of universal quantifiers --- constitute set $\TTT$. Assume these are:
\begin{align*}
Ax_{Ext} :=& \ \forall \yx{0} \yx{1} (\neg \yx{0} = \yx{1} \to \exists \yx{2} \, \yx{2} \eps \yx{0} \triangle \yx{1}).\\
Ax_{Pair} :=& \ \forall \yx{0} \yx{1} \exists \yx{2} \yBpair \text{ with } \yBpair = \forall \yx{3} (\yx{3} \eps \yx{2} \yequi (\yx{3} = \yx{0} \lor \yx{3} = \yx{1})).\\
Empty :=& \ \exists \yx{0} \forall \yx{1} \, \neg \yx{1} \eps \yx{0} \text{ (an instance of the separation schema)}.\\
Ax_{Pow} :=& \ \forall \yx{0} \exists \yx{1} \yBpow \text{ with } \\
& \ \yBpow = \forall \yx{2} ((\yx{2} \eps \yx{1} \to \yx{2} \subseteq \yx{0}) \land (\neg \yx{2} \eps \yx{1} \to \exists \yx{3} \ \yx{3} \eps \yx{2} \setminus \yx{0})).\\
Ax_{Inf} :=& \ \exists \yx{0} (\y{0} \eps \yx{0} \land \forall \yx{1} (\yx{1} \eps \yx{0} \to \exists \yx{2} (\yx{2} \eps \yx{0} \land \yx{2} = \y{suc} \, \yx{1}))).
\end{align*}

Order the formulas in $\hat\TTT$ (consisting of the five axioms and all its subformulas) according to the length of the context in which they occur. Then we find a possible index for which the first case of Definition (\ref{irsetdef}) applies. Note that $\yx{0}$ ranges over $\VVV_{i_0}$, $\yx{1}$ over $\VVV_{i_1}$ and so on. The relevant formulas in $\hat\TTT$ are the existential assertions:
\begin{align*}
0:&\ Empty, \, Ax_{Inf} & i_0=&\,\{0, \omega\} \\
1:&\ \exists \yx{1} \yBpow & i_1=&\,\{1, P\omega\} \\
2:&\ \exists \yx{2} \, \yx{2} \eps \yx{0} \triangle \yx{1}, \, \exists \yx{2} \yBpair, & i_2=&\,\{\, 0, \,\omega, \,2, \{0, P\omega\}, \{\omega, 1\}, \{\omega, P\omega\}\,\}\\
&\ \exists \yx{2} (\yx{2} \eps \yx{0} \land \yx{2} = \y{suc} \, \yx{1}) && \\
3:&\ \exists \yx{3} \, \yx{3} \eps \yx{2} \setminus \yx{0} & i_3=&\,\{\, 0,\, 1,\, \omega\,\}
\end{align*}

$P\omega$ refers to the power set of $\omega$ (which are both single elements in $\VVV$). Let us look more closely at the contexts of length 2: The elements $0$, $\omega$ and $2$ in $i_2$ are witnesses for the formula $\exists \yx{2} \, \yx{2} \eps \yx{0} \triangle \yx{1}$; for instance, $2$ witnesses the difference of $\omega$ and $1$, if $\omega$ is assigned to $\yx{0}$ and $1$ to $\yx{1}$. The pair-sets $2, \{0, P\omega\}, \{\omega, 1\}$ and $\{\omega, P\omega\}$ stem from $\exists \yx{2} \yBpair$, the last formula $\exists \yx{2} \, \yx{2} \eps \yx{0} \land \yx{2} = \y{suc} \, \yx{1}$ only requires the witness $2$ again (due to assigning $\omega$ to $\yx{0}$ and $1$ to $\yx{1}$). These indices satisfy $\nil \ll_\TTT i_0$, $i_0 \ll_\TTT i_1$, $(i_0,i_1) \ll_\TTT i_2$ and $(i_0,i_1,i_2) \ll_\TTT i_3$. The finite model $\VVV_\JJJ$ with $\JJJ = \{i_0,i_1,i_2,i_3,j\}$ and $j = i_0 \cup i_1 \cup i_2 \cup i_3$, is a $\TTT$-submodel of $\VVV_\III$.

Adding more and more axioms and instances of schemata increases the finite model. If we do not add all (actual) infinitely many instances of the schemata at once, we can still use an ``infinite'' element $\omega$ inside the investigated model, but interpret it as potential infinite in the background model by an increasing family of finite, ``real'' sets $\{b \in \VVV_i \mid b \epsi \omega\}$ at stage $i$. For instance, the element $\omega$ at stage $i_2$ is $\{0,2\}$ and $\{0,1 \}$ at stage $i_3$.

\section{Further Remarks}
\label{further}

There are some immediate variations, which we shortly mention.

\subsection{Some Adoptions.}
\label{simpadopt}

An obvious and simple generalization is attained by using a typed or sorted first-order logic. Another generalization is attained by replacing subset inclusion by embeddings $\emb{i}{i'} : \MMM_i \embb\MMM_{i'}$ for $i \leq i'$, which are not necessarily injective. One uses the direct limit instead of the union in all constructions. This includes further examples e.g.~syntactical ones in which terms are identified at a larger stage. 

A transfer to intuitionistic logic with Kripke models is also possible, but requires more effort. This model has two index sets, the preorder of \emph{epistemic states} (the Kripke frame) and the directed set $\III$ of \emph{ontological states}, used here. A relation $R$ in a Kripke structure is then a family of sets $R^k_C$ with $k$ a node and $C$ a context. They satisfy $R^k_C(\vecb{a}) \iff R^{k}_{C'}(\vecb{a})$ for $\vecb{a} \in \MMM_C \cap \MMM_{C'}$ on the one hand and the weaker requirement $R^k_C(\vecb{a})\, \imp\, R^{k'}_C(\vecb{a})$ for $k \leq k'$ and $\vecb{a} \in \MMM_C$ on the other hand.

Notice that we never required the property that a $\ll$-context $(i_0, \dots, i_{n-1})$ satisfies $i_0 \leq \dots \leq i_{n-1}$. This is not necessary, but it can easily be achieved: Add $\dots \cap \bigcap_{0 \leq k < n} \up i_k$, for $C = (i_0,\dots,i_{n-1})$, to the right hand side in the Definition (\ref{muupper}).

One may as well allow terms and functions $f_\HHH := (f_{C \to j})_{Cj \in \HHH}$ with finite maps $f_{C \to j}$ in a straightforward way. The state declaration $\tdt{\y{t}}{C}{j}$ and interpretation $\valp{\y{t}}{j}_{\vecb{a} : C} \in \MMM_j$ for an $n$-ary term $\y{t}$, context $C \in \III^n$, index $j \in \III$ and assignment $\vecb{a} \in \MMM_C$ is
\begin{align*}
\tdt{\y{f}\yt{0} \dots \yt{m-1}}{C}{j} \ &\iffdef \tdt{\yt{0}}{C}{j_0}, \dots, \tdt{\yt{m-1}}{C}{j_{m-1}},\\
\valp{\y{f}\yt{0} \dots \yt{m-1}}{j}_{\vecb{a} : C} \ &:= f_{(j_0,\dots,j_{m-1}) \to j} (\valp{\yt{0}}{j_0}_{\vecb{a} : C}, \dots, \valp{\yt{m-1}}{j_{m-1}}_{\vecb{a} : C}),
\end{align*}
for some $j_0, \dots, j_{m-1}$ with $\tdf{\y{f}}{(j_0,\dots,j_{m-1})}{j} \in \Sigma^\III$. 

Instead of \emph{finite}, we may use the notion of a \emph{definite} collection --- we avoid the notion ``small'' due to its misleading connotation of size. These definite collection are defined by closure properties. An \emph{indefinite} collection is simply a collection that is not definite. Being finite is the least notion of definiteness, whereas the natural reading in set theory is ``definite = size of a set'' and ``indefinite = size of a proper class'', including the example of the cumulative hierarchy $\VVV_{On}$ (with the class of ordinal numbers $On$ as index set and $\VVV_\alpha$ being the $\alpha$'s rank of this hierarchy). A further example is that definite refers to countable sets. In that case we may see the L\"owenheim-Skolem theorem as a special case of the construction described here. For a structure $\MMM_\III$, the index set $\III$ must then be directed with respect to definite sets, that is, for each non-empty definite subset $\JJJ \subseteq \III$ of indices an upper bound of $\JJJ$ exists in $\III$, and $\MMM_\III$ must be \emph{locally definite}, i.e., all sets $\MMM_i$ are definite. With these adoptions, all statements and proofs are carried over easily to this more general situation.

\subsection{Conclusion and Further Work.}
\label{furtherfin}

This paper presents a first step to develop a consequent view of infinity as a potential infinite, that does not require any restrictions of logical inferences. We presented the approach for classical first-order logic as a blueprint for other logics. The core concepts are a formalized notion of an indefinitely large state and a state declaration for expressions. Both concepts allow an interpretation of the universal quantifier that refers to finite sets only.

In order to use it for a larger part of mathematics we have to deal with functions and relations as objects. This requires (potential) infinite objects, which can be accessed only by their approximations. We will extend this approach to a fragment of simple type theory, which includes classical higher-order logic. This requires a more general notion of a system, not only a direct system, and a general notion of a limit of this system. A further challenge is the extension to intuitionistic higher-order logic.

\bibliographystyle{als}
\bibliography{DynModPotInf}

\end{document}